\documentclass{article}

\long\def\rests#1{}

\def\noi{\noindent}

\def\Pf{\noi{\bf Proof.\ \,}}
\def\eop{\hfill\framebox[2.4mm][t1]{\phantom{x}} \vskip 0.15cm } 

\def\voa{vertex operator algebra\ }
\def\gkm{generalized Kac-Moody algebra\ }
\def\gkms{generalized Kac-Moody algebras\ }
\def\fbmla{fake baby monster Lie algebra\ }
\def\a{\alpha}
\def\s{{\bf s}}
\def\g{{\bf g}}

\def\II{I\!I_{1,1}}
\def\VII{V_{\II}}
\def\Z{{\bf Z}}
\def\z{{\bf z}}
\def\wt{{\rm wt}}
\def\stack#1#2{{\scriptstyle #1}\atop{\scriptstyle #2}}
 
\title{A natural construction of Borcherds' Fake Baby Monster Lie Algebra}

\author{Gerald~H\"ohn\thanks{
Mathematisches Institut, Universit\"at Freiburg, 
Eckerstra{\ss}e 1, 79104 Germany. E-mail: {\tt gerald@mathematik.uni-freiburg.de}
}
\ and Nils~R.~Scheithauer\thanks{
Mathematisches Institut, Universit\"at Heidelberg, 
Im Neuenheimer Feld 288, 69120 Heidelberg, Germany. 
E-mail: {\tt nrs@mathi.uni-heidelberg.de}.
\newline
The second author acknowledges financial support from
the Emmy Noether-program of the DFG and
the TMR Network ERB FMRX-CT97-0100 ``Algebraic Lie Theory''.
\newline
\newline}}

\date{25 July 2002}

\begin{document}

\newtheorem{thm}{Theorem}[section]
\newtheorem{prop}[thm]{Proposition}
\newtheorem{lem}[thm]{Lemma}
\newtheorem{rem}[thm]{Remark}
\newtheorem{cor}[thm]{Corollary}
\newtheorem{conj}[thm]{Conjecture}
\newtheorem{de}[thm]{Definition}
\newtheorem{nota}[thm]{Notation}

\maketitle

\begin{abstract}

We use a $\Z_2$-orbifold of the vertex operator algebra associated 
to the Niemeier lattice with root lattice $A_3^8$ and the no-ghost
theorem of string theory to construct a generalized Kac-Moody algebra.
Borcherds' theory of automorphic products allows us to determine the 
simple roots and identify the algebra with the fake baby monster Lie algebra.

\end{abstract}


\section{Introduction}

Up to now, there are only three generalized Kac-Moody algebras or superalgebras
for which natural constructions are known.
These are the fake monster Lie algebra~\cite{B-fake} and the monster
Lie algebra~\cite{B-moonshine} constructed by Borcherds and the fake monster Lie 
superalgebra~\cite{S-superfake} constructed by the second author.
All these algebras can be interpreted as the 
physical states of a string moving on a certain target space.

In~\cite{B-moonshine}, there was also introduced a method to obtain
new \gkms from old ones by twisting the denominator identity with some
outer automorphism.
Such Lie algebras are only defined through generators and relations,
as it is the case for all other known examples of \gkms
(see, for e.g.,~\cite{GN}). 
In particular, Borcherds found a \gkm of rank~$18$ called the {\em \fbmla} 
by taking a ${\bf Z}_2$-twist of the fake monster Lie algebra 
(see~\cite{B-moonshine}, Sect.~14, Example~1)
and he asked for a natural construction of it. 
The purpose of this note is to present such a construction.

\medskip

The fake monster and the monster Lie algebra are obtained 
in the following way: Take for $V$ the \voa (VOA) $V_{\Lambda}$ 
associated to the Leech lattice $\Lambda$ or the Moonshine module VOA 
$V^{\natural}$ and let $\VII$ be the vertex algebra of the
two dimensional even unimodular Lorentzian lattice $\II$. The tensor product
$V\otimes \VII$ is a vertex algebra of central charge $26$ with
an invariant nonsingular bilinear form.
Let $P_n$ be the subspace of Virasoro highest weight vectors of 
conformal weight $n$, i.e., the space of vectors $v$ satisfying
$L_m(v)=0$ for $m>0$ and $L_0(v)=n\cdot v$. Then $P_1/L_{-1}P_0$ is
a Lie algebra with an induced invariant bilinear form $(\,,\,)$. 
The fake monster or monster Lie algebra is defined as the quotient
of $P_1/L_{-1}P_0$ by the radical of $(\,,\,)$.
The non-degeneracy of the induced bilinear form is used
to show that one has indeed obtained a generalized Kac-Moody algebra $\g$.

Alternatively, one can use the bosonic ghost vertex superalgebra 
$V_{\rm ghost}$ of central charge $-26$ and define the
Lie algebra $\g$ as the BRST-cohomology group 
$H^1_{\rm BRST}(V\otimes \VII)$ (cf.~\cite{FGZ}).

\smallskip

The above construction can be carried out for any VOA $V$ of central
charge~$24$, but until now only for $V_{\Lambda}$ and $V^{\natural}$
it was known how to compute the simple roots of the \gkm $\g$.
The Lie algebras obtained
from the lattice VOAs $V_K$, where $K$ is any rank $24$ Niemeier lattice,
are all isomorphic to the fake monster Lie algebra, because 
$K\oplus \II$ is always equal to the even unimodular Lorentzian 
lattice $I\!I_{25,1}$.
The Moonshine module $V^{\natural}$ was constructed in~\cite{FLM}
as a $\Z_2$-orbifold of $V_{\Lambda}$. This $\Z_2$-orbifold construction was
generalized to any even unimodular lattice instead of $\Lambda$
in~\cite{DGM,DGM2}.

In our construction of the fake baby monster Lie algebra
we take for $V$ the $\Z_2$-orbifold of $V_K$,
where $K$ is the Niemeier lattice with root lattice~$A_3^8$.
The computation of the root multiplicities is harder than in the
previous cases:
The weight $1$ part $V_1$ of $V$ is the semisimple Lie algebra of type
$A_1^{16}$ and~$V$ forms an integrable highest weight representation of
level~$2$ for the affine Kac-Moody algebra of type $\widehat{A}_1^{16}$.
The decomposition of $V$ into $\widehat{A}_1^{16}$-modules can be 
described by the Hamming code ${\cal H}_{16}$ of length~$16$ and its
dual. We use this combinatorial description together with the no-ghost
theorem to determine the root lattice and root multiplicities
of $\g$. The multiplicities 
obtained are exactly the exponents of a product expansion
of an automorphic form constructed in~\cite{S-twist}.
This allows us to interpret the automorphic product as 
one side of the denominator identity of $\g$,
to determine its simple roots, and finally
to identify $\g$ with the fake baby monster Lie algebra.

\medskip

The paper is organized in the following way.  
In Section~\ref{voav}, the construction of the \voa $V$ is described. 
We use a formula of Kac and Peterson (cf.~\cite{Kac}, Ch.~13) to express the
character of the affine Kac-Moody algebra of type $\widehat{A}_1^{16}$
through string functions and theta series.
In the last section,
the root lattice and root multiplicities of $\g$ are computed 
and $\g$ is identified as the fake baby monster Lie algebra.

\medskip

We would like to thank R.~Borcherds for helpful comments on an 
early version of this paper.


\section{The VOA $V$}\label{voav}

We define a VOA $V$ of central charge~$24$ and compute its
character as representation for an affine Kac-Moody algebra.

\smallskip

Recall that there exist exactly $24$ positive 
definite even unimodular lattices in dimension~$24$. They can be
classified by their root sublattice.
\begin{de}\rm
Let $V=V_K^+\oplus (V_K^T)^+$ be the $\Z_2$-orbifold of the lattice
VOA associated to the Niemeier lattice $K$ with root lattice $A_3^8$.
Here, $T$ is the involution in  ${\rm Aut}(V_K)$ which is the up to 
conjugation unique lift of the involution $-1$ in ${\rm Aut}(K)$ 
to ${\rm Aut}(V_K)$ (cf.~\cite{DGH}, Appendix~D);  $V_K^+$ is the
fixpoint subVOA of $V_K$ under the action of $T$ 
and $(V_K^T)^+$ is the fixpoint set of the $T$-twisted module~$V_K^T$. 
\end{de}

In~\cite{DGM,DGM2} it is shown that $V$ has the structure
of a VOA of central charge~$24$.

\begin{thm}
Let $V_{A_{1,2}}$ be the VOA which has the integrable
level~$2$ representation of highest weight $(2,0)$
for the affine Kac-Moody algebra of type $\widehat{A}_1$ as
underlying vector space.
The subVOA $\widetilde{V}_{1}$ generated by the weight $1$ subspace
$V_1$ of $V$ is isomorphic to the affine Kac-Moody VOA $V_{A_{1,2}^{16}}$, 
the tensor product of $16$ copies of $V_{A_{1,2}}$.
\end{thm}

\Pf \cite{DGM}. \eop
\begin{rem}\rm
As noted in the introduction of~\cite{GH}, $V$ is the unique
VOA in the genus of the Moonshine module containing 
$V_{{A}_{1,2}^{16}}$ as a subVOA. (See~\cite{H} for the definition
of the genus of a VOA.) 
\end{rem}

\smallskip

To describe the decomposition of $V$ as a $V_{{A}_{1,2}^{16}}$-module 
in a convenient way and for some later applications, we explain
some well known properties of the binary Hamming code ${\cal H}_{16}$
of length~$16$.

Let ${\cal H}_{16}^{\perp}\subset {\bf F}_2^{16}$ be the binary code
spanned by the rows of the matrix
$$\left(\begin{array}{c}
1\,  1\,  1\,  1\,  1\,  1\,  1\,  1\,  1\,  1\,  1\,  1\,  1\,  1\,  1\,  1\, \\
1\,  1\,  1\,  1\,  1\,  1\,  1\,  1\,  0\,  0\,  0\,  0\,  0\,  0\,  0\,  0\, \\
1\,  1\,  1\,  1\,  0\,  0\,  0\,  0\,  1\,  1\,  1\,  1\,  0\,  0\,  0\,  0\, \\
1\,  1\,  0\,  0\,  1\,  1\,  0\,  0\,  1\,  1\,  0\,  0\,  1\,  1\,  0\,  0\, \\
1\,  0\,  1\,  0\,  1\,  0\,  1\,  0\,  1\,  0\,  1\,  0\,  1\,  0\,  1\,  0\, 
\end{array}\right)_.$$ 
This code is known as the first order Reed-Muller code of length~$16$
and has ${\rm AGL}(4,2)$ as automorphism group (cf.~\cite{DGH}, Th.~C.3).
Its Hamming weight enumerator is 
$W_{{\cal H}_{16}^{\perp}}(x,y)=x^{16}+ 30\,x^8y^8+ y^{16}$.
The Hamming code ${\cal H}_{16}$ is defined as the dual code of
${\cal H}_{16}^{\perp}$: these are the vectors $c\in {\bf F}_2^{16}$ satisfying
\hbox{$\sum_{i=1}^{16}c_i\cdot d_i=0$} for all $d\in {\cal H}_{16}^{\perp}$. 
We easily see that ${\cal H}_{16}$ can correct $1$-bit errors and so the smallest
nonzero code vector has weight at least~$4$. Indeed, by
the MacWilliams identity we obtain for its weight enumerator:
$$ W_{{\cal H}_{16}}(x,y)=
x^{16}+y^{16}+140\,(x^4y^{12}+x^{12}y^4)+448\,(x^6y^{10}+x^{10}y^6)+870\,x^8y^8.$$
It follows that the $140$ codewords of weight~$4$ form a Steiner system
of type $S(16,4,3)$, i.e., for every $3$-tuple of coordinate positions there is
exactly one weight~$4$ code vector with value $1$ at this $3$ positions.

We also need the weight enumerator of all other cosets in the cocode 
${\bf F}_2^{16}/{\cal H}_{16}$. The $2^5$ cosets ${\cal H}_{16}+c$ 
can be represented by vectors $c$ of type $(0^{16})$ (one coset), 
$(0^{15}1^1)$ (sixteen cosets) and $(0^{14}1^2)$ (fifteen cosets). Indeed,
the vectors of type $(0^{16})$ and $(0^{15}1^1)$ must be in different cosets
and for every vector of type $(0^{14}1^2)$ there are by the Steiner system
property exactly $7$ others contained in the same coset. Since 
$1+16+{16 \choose 2}/8 =2^5$, all cosets are counted.
For the cosets of type $(0^{15}1^1)$, the weight enumerator is
$$W_{{\cal H}_{16}+c}(x,y)=\frac{1}{32}\big((x+y)^{16}-(x-y)^{16}\big)$$
$$ = x^{15}y + xy^{15} + 35\, (x^{13}y^3+x^3y^{13} )+ 
   273\, (x^{11}y^5+x^5y^{11}) + 715\,(x^9y^7+x^7y^9) $$
because these cosets contain only vectors of odd weight and ${\rm AGL}(4,2)$
acts transitively on the coordinates.
For the cosets of type $(0^{14}1^2)$, the weight enumerator is
$$W_{{\cal H}_{16}+c}(x,y)=
\frac{1}{15}\big(((x+y)^{16}+(x-y)^{16})/2-W_{{\cal H}_{16}}(x,y)\big)$$
$$=8\,(x^{14}y^2 +x^2y^{14}) + 112\,(x^{12}y^4+x^4y^{12}) + 
    504\,(x^{10}y^6+x^6y^{10}) + 800\,x^8y^8 $$
because ${\rm AGL}(4,2)$ acts also transitively on pairs of coordinates.

\smallskip

We need two other weight enumerators. 
Let $({\bf F}_2^8)_{0}$ be the subcode 
of all vectors of even weight in ${\bf F}_2^8$ 
and $({\bf F}_2^8)_{1}$ be the coset of vectors of odd weight.
Their weight enumerators are:
\begin{eqnarray*}
W_{({\bf F}_2^8)_{0}}(x,y) \!\! & \!=\! &  \! \!\frac{1}{2}\big((x+y)^8+(x-y)^8\big)
=x^8+y^8+28\,(x^6y^2+x^2y^6)+70\,x^4y^4, \\
W_{({\bf F}_2^8)_{1}}(x,y)\!\! &\! =\! & \!\! \frac{1}{2}\big((x+y)^8-(x-y)^8\big)
= 8\,(x^7y+y^7x) + 56\,(x^5y^3+x^3y^5).
\end{eqnarray*} 

\medskip 

The rational Kac-Moody VOA $V_{A_{1,2}}$ has three irreducible modules 
$M(0)$, $M(1)$ and $M(2)$ of conformal weight $0$, $1/2$ and 
$3/16$, respectively. The irreducible modules of
$V_{A_{1,2}^{16}}\cong V_{A_{1,2}}^{\otimes 16}$ are the tensor products
$M({i_1})\otimes\cdots\otimes
M({i_{16}})$, $i_1$, $\ldots$, $i_{16}\in\{0,\,1,\,2\}$, for which we
write shortly $M(i_1,\ldots\,,i_{16})$.

\begin{thm}\label{Vdecomp}
Up to permutation of the $16$ tensor factors $V_{A_{1,2}}$,
the $V_{A_{1,2}^{16}}$-module decomposition of $V$ 
into isotypical components has the following structure: 
$$V=\bigoplus_{\delta\in {\cal H}_{16}^{\perp}}K(\delta),$$
where $K(\delta)$ is defined 

\begin{tabular}{ll}
for  $\wt(\delta)=0$ by &
$\bigoplus\limits_{c\in {\cal H}_{16}}M(c)$,   \\
for  $\wt(\delta)=8$ by &
$\bigoplus\limits_{i_1,\ldots,i_{16}\in\{0,1,2\} } 
n_{i_1,\ldots,i_{16}}^{\delta} \, M(i_1,\ldots,i_{16})$,     \\
\noalign{\noindent
here, $n_{i_1,\ldots,i_{16}}^{\delta}=1$, if $[i_k/2]=\delta_k$
(where $[x]$ denotes the Gau\ss  bracket of $x$) 
\newline 
for $k=1$, $\ldots$, $16$ and $\# \{k \mid i_k=1\}$ is odd,
and $n_{i_1,\ldots,i_{16}}^{\delta}=0$ otherwise,\phantom{$\frac{|}{|}$}  
 } 
and for $\wt(\delta)=16$ by  &
$2^3\, M(2,\ldots,2)$. \phantom{$\frac{|}{|}$}
\end{tabular}

\end{thm}

\Pf 
This follows by applying a variation of Th.~4.7 in~\cite{DGH}
to the glue code $\Delta$ of the Niemeier lattice with root lattice $A_3^8$.
The Virasoro VOA $L_{1/2}$ of central charge $1/2$ there
is replaced by the VOA $V_{A_{1,2}}$ which has an 
isomorphic fusion algebra and the lattice $D_1$ is replaced by the 
lattice $A_3$ which has an isomorphic discriminant group ${\bf Z}/4{\bf Z}$.
Then the theorem remains valid if one substitutes the three irreducible
$L_{1/2}$-modules of weight $0$, $1/2$ and $1/16$ by the three 
irreducible $V_{A_{1,2}}$-modules of weight $0$, $1/2$ and $3/16$, respectively.
The explicit description of the decomposition resulting from the 
above ${\bf Z}/4{\bf Z}$-code~$\Delta$ of length~$8$ 
was given in the proof of Th.~5.3 in~\cite{DGH}. 
See also the following Remark~5.4 (2) there. \phantom{xxxx}
\eop

\medskip

The ${\bf Z}$-grading on the VOA $V=\bigoplus_{n=0}^{\infty}V_n$ is given 
by the eigenvalues of the Virasoro generator~$L_0$.
There is also the action of the Lie algebra of type $A_1^{16}$.
For $\s$ in the weight 
lattice~$(A'_1)^{16}\cong\big(\frac{1}{\sqrt{2}}{\bf Z}\big)^{16}$
we denote by $V_n(\s)$ the subspace of $V_n$ on which the action of
the Cartan subalgebra of the Lie algebra $A_1^{16}$
has weight $\s$. The character of~$V$ defined by 
$$\chi_V=q^{-1}\sum_{n\in \bf Z}\sum_{\s\in (A'_1)^{16}}
          \dim V_n(\s)\,q^n\,e^{\s}$$
is an element in the ring of formal Laurent series in $q$ with
coefficients in the group ring ${\bf C}[(A'_1)^{16}]$.

For the proof of some identities, it is useful to
interpret an element $f$ in ${\bf C}[L][[q^{1/k}]][q^{-1/k}]$,
where $L$ is a lattice and $k\in{\bf N}$, 
as a function on ${\cal H}\times (L\otimes {\bf C})$,
where ${\cal H}=\{z\in {\bf C} \mid \Im(z)>0\}$ is the complex
upper half plane. This is done by the substitutions 
$q\mapsto e^{2\pi i \tau}$ and $e^{\s} \mapsto e^{2\pi i(\s,\z)}$
for $(\tau,\z)\in {\cal H}\times (L\otimes {\bf C})$ 
(in the case of convergence). We indicate this by writing $f(\tau,\z)$.

\smallskip

To compute $\chi_V$ with the help of the Weyl-Kac character formula, 
we need various power series, which are the Fourier expansion of 
various modular and Jacobi forms.
Let $\eta(\tau)=q^{1/24}\prod_{n=1}^{\infty}(1-q^n)$ 
be the Dedekind eta function.

First, there are the three ``string functions"´ $c_0$, $c_1$ and $c_2$ 
which are modular functions for $\widetilde\Gamma(16)$ of weight $-1/2$.  
They are defined by
\begin{eqnarray*}
c_0 \! &\! =\! & \!\frac{1}{2}\left( \frac{\eta(\tau/2)}{\eta(\tau)^2}
+\frac{\eta(\tau)}{\eta(2\tau)\eta(\tau/2)}\right)
=\, q^{-1/16}\cdot(
1 + q + 3\,q^2 + 5\,q^3 + 10\,q^4 + \cdots
),\\
 c_1\! &\! =\! &\! \frac{1}{2}\left( \frac{\eta(\tau/2)}{\eta(\tau)^2}
-\frac{\eta(\tau)}{\eta(2\tau)\eta(\tau/2)}\right)
=\,q^{-1/16}\cdot(
 q^{1/2} + 2\,q^{3/2} + 4\,q^{5/2} + \cdots
), \\
c_2\! &\! =\! &\! \frac{\eta(2\tau)}{\eta(\tau)^2}=\, 1 + 2\,
  q + 4\,{q}^2 + 8\,{q}^3 + 
   14\,{q}^4 + 24\,{q}^5 + 40\,{q}^6 + \cdots\,.
\end{eqnarray*}

We denote the dual lattice of an integral lattice $L$ by $L'$.
As we are dealing with level~$2$ representations of $\widehat{A}_1$,
it will be convenient to define for $\gamma\in L'/L$
the theta function by
$$\Theta_{L+\gamma}=
       \sum_{{\bf s}\in L+\gamma} q^{\s^2/4}\,e^{\bf s},$$
writing $\s^2=(\s,\s)$ for the norm of $\s$ and
$L+\gamma$ for the coset of $L$ in $L'$ determined by $\gamma$.
In particular, we define the following three theta series: 
\begin{eqnarray*}
\vartheta_0(\tau,z) & = & \Theta_{2A_1}(\tau,z),\\
\vartheta_1(\tau,z) & = &  \Theta_{2A_1+ \sqrt{2}}(\tau,z),\\
\vartheta_2(\tau,z) & = &  \Theta_{A_1+ \frac{1}{\sqrt{2}}}(\tau,z)
\ = \  \Theta_{2A_1+ \frac{1}{\sqrt{2}}}(\tau,z) +
\Theta_{2A_1- \frac{1}{\sqrt{2}}}(\tau,z).
\end{eqnarray*}

The graded characters  $\chi_i=\chi_{M(i)}$
of the three irreducible level $2$ representations 
$M(i)$, $i=0$, $1$ and $2$,
of the affine Kac-Moody Lie algebra of type~$\widehat{A}_1$ can now be
expressed in the above series:
\begin{prop}[Kac-Peterson{\rm, cf.~\cite{Kac}, Ch.~13}]\label{strings}
\begin{eqnarray*}
\chi_0 & = & \chi_{M(0)} = c_0\cdot\vartheta_0+c_1\cdot\vartheta_1, \\
\chi_1 & = & \chi_{M(1)} = c_1\cdot\vartheta_0+c_0\cdot\vartheta_1, \\
\chi_2 & = & \chi_{M(2)} = c_2\cdot\vartheta_2.
\end{eqnarray*} \eop
\end{prop}

We combine this information with the $V_{A_{1,2}^{16}}$-decomposition
of $V$.

\begin{prop}\label{kac-moody-character-codes}
For $\delta \in {\bf F}_2^{16}$ and $d\in {\bf F}_2^{n-\wt(\delta)}$
(identified with the subspace 
$\{c\in  {\bf F}_2^{16} \mid c_i=0\ \hbox{for all\ } i=1,\ldots,16
\hbox{\ with\ }\delta_i=1\}$)
we introduce the shorthand notation
$$\vartheta_d^{\delta}(\tau,\z)=
\prod_{\stack{i\in\{1,\ldots,16\}}{\delta_i=0}}\vartheta_{d_i}(\tau,z_i)
\quad\hbox{and}\quad
\vartheta_2^{\delta}(\tau,\z) =
\prod_{\stack{i\in\{1,\ldots,16\}}{\delta_i=1}}\vartheta_{2}(\tau,z_i).$$
Then the character of $V$ is
\begin{eqnarray*}
\chi_V(\tau,\z) & = &
\sum_{d\in {\bf F}_2^{16}}   
   W_{{\cal H}_{16}+d}(c_1,c_2)\,\vartheta_d^{(0,\,\ldots,\,0)}(\tau,\z)\\
&& + \sum_{\stack{\delta \in  {\cal H}_{16}^{\perp}}{ {\rm wt}(\delta)=8}}
 \sum_{d\in {\bf F}_2^8}
 W_{({\bf F}_2^{8})_1+d}(c_1,c_2)\,\vartheta_d^{\delta}(\tau,\z)
\,\cdot \, c_2^8\, \vartheta_2^{\delta}(\tau,\z) \\
&& + 2^3c_2^{16}\, \vartheta_2^{(1,\,\ldots,\,1)}(\tau,\z).
\end{eqnarray*}
\end{prop}

\Pf  This is a rather trivial resummation. 
{}From Theorem~\ref{Vdecomp} and the notation there one gets
$$\chi_V(\tau,\z) = 
 \sum_{c\in {\cal H}_{16}} \chi_{M(c)}
+\sum_{\stack{\delta\in{\cal H}_{16}^{\perp}}{\scriptstyle \wt(\delta)=8}}
    \sum_{\stack{ i_1,\ldots,i_{16}}{\phantom{\frac{1}{1}}\in\{0,1,2\}}}
n_{i_1,\ldots,i_{16}}^{\delta}\chi_{M(i_1,\ldots,i_{16})}
 +2^3\chi_{M(2,\ldots,2)}$$
by Proposition~\ref{strings}
\begin{eqnarray*}
&  = &\sum_{c\in {\cal H}_{16}}\sum_{d\in {\bf F}_2^{16}}
c_0^{16-\wt(c+d)}c_1^{\wt(c+d)}\prod_{i=1}^{16}\vartheta_{d_i}(\tau,z_i)\\
& &  +\sum_{\stack{\delta\in{\cal H}_{16}^{\perp}}{ \wt(\delta)=8}}
\sum_{c \in ({\bf F}_2^8)_1}\sum_{d \in {\bf F}_2^8}
c_0^{8-\wt(c+d)}c_1^{\wt(c+d)}
\ \Big(\prod_{\stack{i\in\{1,\ldots,16\}}{\delta_i=0}}
\vartheta_{d_i}(\tau,z_i)\Big)\\
& &  \qquad\qquad\qquad\qquad\qquad\qquad\qquad\qquad\qquad\qquad
\cdot\ c_2^8\ \Big( \prod_{\stack{i\in\{1,\ldots,16\}}{\delta_i=1}}
\vartheta_{2}(\tau,z_i)\Big) \\
& &\  +2^3c_2^{16}\prod_{i=1}^{16}\vartheta_{2}(\tau,z_i),
\end{eqnarray*}
which simplifies to the formula given in the proposition.
\eop

We can simplify this expression further. 
We define three modular functions $h$, $g_0$ and $g_1$ of weight $-8$ for 
the modular group $\Gamma(2)$:
\begin{eqnarray*}
h(\tau)\! &\! =\! &\! \eta(\tau)^{-8}\eta(2\tau)^{-8}=
q^{-1}+  8+ 52\,q+ 256\, q^2 + 1122\, q^3 + 4352\, q^4 + \cdots, \\
g_0(\tau)\!&\! =\! &\! \frac{1}{2}\big(h(\tau/2)+h((\tau+1)/2)\big)=
8 + 256\, q + 4352\, q^2 +  52224\,q^3 +\cdots,  \\
g_1(\tau)\! &\! =\! &\!\frac{1}{2}\big(h(\tau/2)-h((\tau+1)/2)\big)=
 q^{-1/2}+ 52\, q^{1/2} + 1122\,  q^{3/2} + \cdots.
\end{eqnarray*}

\begin{lem}\label{modulids}
One has the following identities between modular functions
for~$\widetilde\Gamma(16)$:
\begin{eqnarray*}
g_0+h &= & W_{{\cal H}_{16}}(c_1,c_2), \\
g_0 & = & W_{{\cal H}_{16}+(1,1,0,\ldots,0)}(c_1,c_2) 
\ =\ W_{({\bf F}_2^{8})_1}(c_1,c_2)\cdot c_2^8
\ =\ 2^3\,c_2^{16},\\
g_1 & = & W_{{\cal H}_{16}+(1,0,0,\ldots,0)}(c_1,c_2)
\ =\ W_{({\bf F}_2^{8})_0}(c_1,c_2)\cdot c_2^8.
\end{eqnarray*}
\end{lem}

\Pf The space of modular functions for $\widetilde\Gamma(16)$ 
of weight $-8$ with poles of given order only at the cusps is 
finite dimensional. Comparing the Fourier expansions on both 
sides gives the result. The details are left to the reader.
\eop

Let $\pi: A_1^{16} \longrightarrow A_1^{16}/ (2A_1)^{16} \cong {\bf F}_2^{16}$
be the projection map. We define the even integral lattice 
$N=\frac{1}{\sqrt{2}}\pi^{-1}({\cal H}_{16})$.
Its discriminant group $N'/N$ is contained in
$(\frac{1}{\sqrt{2}} A_1')^{16}/(\frac{1}{\sqrt{2}}2A_1)^{16}
\cong {\bf Z}_4^{16}$ and fits into the short exact sequence
$$0\longrightarrow {\bf F}_2^{16}/{\cal H}_{16}
\stackrel{\iota}{\longrightarrow} 
N'/N \longrightarrow {\cal H}_{16}^{\perp}\longrightarrow 0.$$
Since ${\cal H}_{16}^{\perp}\subset {\cal H}_{16}$ one has
$N'/N\cong {\bf Z}_2^{10}$, the sequence has a (non-canonical) split 
$\mu: {\cal H}_{16}^{\perp}\longrightarrow N'/N$
and all the squared lengths of $N'$
are integral, i.e., the induced discriminant form is of type $2^{+10}_{I\!I}$. 

For $\gamma\in N'/N$ define the function
$$ f_{\gamma}=\cases{
     g_0 +h, & if $\gamma=0$, \cr
     g_0, & if \, $\gamma^2\equiv 0 \pmod{2}$ and $\gamma \not=0$,\cr
     g_1, & if \,  $\gamma^2\equiv 1 \pmod{2}$.}
\qquad\qquad(*)   $$

We collect now all the theta functions $\vartheta_d^{\delta}(\tau,{\bf z})$ 
and $\vartheta_2^{\delta}(\tau,{\bf z})$ in $\chi_V$
and arrive at our final expression for the character of $V$.

\begin{thm}\label{kac-moody-character-lattice}
$$ \chi_V(\tau,{\bf z}) = 
  \sum_{\gamma\in N'/N} 
f_{\gamma}(\tau)\, \Theta_{\sqrt{2}(N+\gamma)}(\tau,{\bf z}). $$
\end{thm}
\Pf The expression for 
$\chi_V(\tau,{\bf z})$ obtained in Proposition~\ref{kac-moody-character-codes}
can be rewritten as
\begin{eqnarray*}
\chi_V(\tau,{\bf z}) & = & 
\sum_{d\in {\bf F}_2^{16}}  W_{{\cal H}_{16}+d}(c_1,c_2)\,\cdot \,
\Theta_{\sqrt{2}(\pi^{-1}(d))}(\tau,{\bf z}) \\
&& + \sum_{\stack{ \delta \in  {\cal H}_{16}^{\perp}}{ {\rm wt}(\delta)=8}}
 \sum_{d\in {\bf F}_2^{16}}
 W_{({\bf F}_2^{8})_1+d'}(c_1,c_2)\,\cdot \, c_2^8\,\,\cdot \,
\Theta_{\sqrt{2}(\pi^{-1}(d)+\mu(\delta))}(\tau,{\bf z}) \\
&& + \sum_{d\in {\bf F}_2^{16}}\  2^3c_2^{16}\,\cdot \,
\Theta_{\sqrt{2}(\pi^{-1}(d)+\mu(1,\,\ldots,\,1))}(\tau,{\bf z}), \\
\noalign{
\noindent where $d'$ is the vector of those components $d_i$ of $d$ 
for which $\delta_i=0$. 
Using Lemma~\ref{modulids}, one gets}\\
 & = & 
\sum_{\bar d\in {\bf F}_2^{16}/{\cal H}_{16}} f_{\iota(\bar d)}(\tau)
\,\cdot \, \Theta_{\sqrt{2}(N+\iota(\bar d))}(\tau,{\bf z}) \\
&& + \sum_{\stack{ \delta \in  {\cal H}_{16}^{\perp}}{ {\rm wt}(\delta)=8}}
 \sum_{\bar d\in {\bf F}_2^{16}/{\cal H}_{16}}
f_{\iota(\bar d)}(\tau) \,\cdot \,
\Theta_{\sqrt{2}(N+\iota(\bar d)+\mu(\delta))}(\tau,{\bf z}) \\ \nopagebreak[4] 
&&+\sum_{\bar d\in {\bf F}_2^{16}/{\cal H}_{16}}
 f_{\iota(\bar d)}(\tau)\,\cdot \,
\Theta_{\sqrt{2}(N+\iota(\bar d) +  \mu(1,\,\ldots,\,1))}(\tau,{\bf z}) \\ 
\pagebreak[4]
& = &  \sum_{ \delta \in  {\cal H}_{16}^{\perp}}
\sum_{\bar d\in {\bf F}_2^{16}/{\cal H}_{16}}
f_{\iota(\bar d)}(\tau)\,\cdot \,
\Theta_{\sqrt{2}(N+\iota(\bar d)+\mu(\delta))}(\tau,{\bf z}) \\
& = &  \sum_{\gamma\in N'/N} 
f_{\gamma}(\tau)\,\cdot\, \Theta_{\sqrt{2}(N+\gamma)}(\tau,{\bf z}).
\qquad\qquad\qquad\qquad\qquad\qquad
\hbox{\framebox[2.4mm][t1]{\phantom{x}}} 
\end{eqnarray*}


\section{The \gkm $\g$}\label{fbmla}

In this section,
we use the no-ghost theorem to construct a \gkm $\g$ from $V$.
Theorem~\ref{kac-moody-character-lattice} allows us to describe its root 
multiplicities. We determine the simple roots
with the help of the singular theta-correspondence and show that $\g$ is
isomorphic to the fake baby monster Lie algebra.

\medskip

There is an action of the BRST-operator on the tensor product 
of a vertex algebra $W$ of central charge~$26$
with the bosonic ghost vertex superalgebra
$V_{\rm ghost}$ of central charge~$-26$,  
which defines the BRST-cohomology groups  $H^*_{\rm BRST}(W)$. 
The degree $1$ cohomology group $H^1_{\rm BRST}(W)$ has additionally 
the structure of a Lie algebra, see~\cite{FGZ,LZ,Z}. 

Let $V$ be the VOA of the last section. 
As it is the case for the Moonshine module, we can assume that $V$ is 
defined over the field of real numbers. 
The same holds for the vertex algebra  $\VII$  associated to the even unimodular 
Lorentzian lattice $\II$ in dimension~$2$ and for $V_{\rm ghost}$.
\begin{de}\rm
We define the Lie algebra $\g$ as  $H^1_{\rm BRST}(V\otimes \VII).$
\end{de}

Let $L=N\oplus \II$, where $N$ is the even lattice 
defined in the previous section. 

\begin{prop}\label{rootspace}
The Lie algebra $\g$ is a \gkm graded by the lattice $N'\oplus \II=L'$.
Its components $\g(\a)$, for $\a=(\s,r)\in N'\oplus \II$ are
isomorphic to $V_{1-r^2/2}(\sqrt{2}\,\s)$ for $\a\not = 0$ and to 
$V_1(0)\oplus {\bf R}^{1,1}\cong {\bf R}^{17,1}$ for $\a=0$.
\end{prop}

\Pf 
The vertex algebra $V \otimes V_{\II}$ has a canonical 
invariant bilinear form which can be used to show that the construction 
of $\g$ as BRST-cohomology  group is equivalent to the so called old covariant 
construction used in~\cite{B-moonshine} 
(cf.~\cite{LZ}, section~2.4, cf.~\cite{Z}, section~4).  
In more detail, $V$ carries an action of the Virasoro algebra of central 
charge~$24$ and has positive definite bilinear form such 
that the adjoint of the Virasoro generator $L_n$ is $L_{-n}$ (see~\cite{DGM}). 
Similarly, $V_{\II}$ has an invariant bilinear form 
(cf.~\cite{S98}, section~2.4), and on $V \otimes 
V_{\II}$, we take the one induced from the tensor product.
This allows us to work in the old covariant picture.

The second part now follows from the no-ghost theorem as given 
in~\cite{B-moonshine}, Th.~5.1, if we use for $G$ a maximal torus of the
real Lie group $SU(2)^{16}$ acting on $V$. 
The proof of the first part is similar to that of Th.~6.2. of~\cite{B-moonshine}.
\eop

The subspace $\g(0)$ of degree $0 \in L'$ is a Cartan subalgebra for $\g$.

\begin{thm}\label{rootmult}
Let $f_{\gamma}(\tau)=\sum_{n\in\Z}a_{\gamma}(n)\,q^n$ be the 
Fourier expansion of the $f_{\gamma}$, \hbox{$\gamma \in N'/N$},
defined in~$(*)$.
For a nonzero vector $\a\in L'$
the dimension of the component $\g(\a)$ is given by
$$ \dim \g(\a)= a_{\gamma}(-\alpha^2/2),   $$
where $\gamma$ is the rest class of $\a$ in $L'/L\cong N'/N$. 
The dimension of the Cartan subalgebra is $18$.
\end{thm}

\Pf Theorem~\ref{kac-moody-character-lattice}, Proposition~\ref{rootspace}.
\eop 

It follows from the Fourier expansion of $f_\gamma$
that the real roots of $\g$ are the norm~$1$ vectors in $L'$ and 
the norm~$2$ vectors in $L$ both with multiplicity~$1$. The real roots of  $\g$ 
generate the Weyl group $W$ of $\g$ which is also equal to the reflection 
group of $L'$. Hence the real simple roots of $\g$ are the simple roots of 
the reflection group of $L'$. 

\begin{prop}
There is a primitive norm $0$ vector $\rho$ in $L'$, called the Weyl vector, 
such that the simple roots of the reflection group of $L'$ are
the roots $\a$ satisfying $(\rho,\a)=-\a^2/2$.
\end{prop}

\Pf Let $\Lambda_{16}$ be the Barnes-Wall lattice.
We write $L(k)$ for the lattice obtained from the lattice $L$
by rescaling all norms by a factor $k$.
Since the discriminant forms 
of the lattices $L=N\oplus\II$ and $\Lambda_{16}\oplus \II(2)$ 
are equal, both lattices are in the genus $I\!I_{17,1}(2^{+10}_{I\!I})$. 
It follows from Eichler's theory of spinor genera that there is 
only one class in this genus and so both lattices must be isomorphic. 
For the rescaled dual of the Barnes-Wall lattice we have 
$\Lambda_{16}'(2)\cong\Lambda_{16}$ 
so that $L'(2)\cong\Lambda_{16}\oplus \II$.
The reflection  group of $\Lambda_{16}\oplus \II$ has a primitive norm~$0$ 
vector $\rho$ such that the simple roots are the roots satisfying $(\rho,\a)=-\a^2/2$ 
(e.g.~\cite{B-theta}, Example~12.4). This implies the statement. 
\eop 

\begin{rem}\rm
\begin{enumerate}
\item If we write $L=\Lambda_{16}\oplus \II(2)$ with elements $(\s,m,n)$, 
\hbox{$s\in \Lambda_{16}$}, $m$, $n\in \Z$ and norm
$(\s,m,n)^2=\s^2-4mn$ we can take $\rho=({\bf 0},0,1/2)$.
Then the simple roots of the reflection group of $L'$ are the 
norm~$1$ vectors in $L'$ of the form $(\s,1/2,(s^2-1)/2)$, $\s\in \Lambda_{16}'$,
and the norm~$2$ vectors $(\s,1,(\s^2-2)/4)$ in $L$, i.e., $\s\in\Lambda_{16}$
with $4|(\s^2-2)$.
\item The automorphism group $\mbox{Aut}(L')^+$ is the semidirect product of the 
reflection subgroup by a group of diagram automorphisms. 
Since $\Lambda_{16}$ has no roots, Theorem~3.3 of~\cite{B-leechlike} 
implies that the group of diagram automorphisms is equal 
to the group of affine automorphisms of the Barnes-Wall lattice. See also ~\cite{B-lorentz}, p.~345.
\end{enumerate}
\end{rem}
We fix a Weyl vector $\rho$ and the Weyl chamber containing $\rho$. 
\begin{prop}
The positive multiples $n\rho$ of the Weyl vector are imaginary simple roots 
of $\g$ with multiplicity $16$ if $n$ is even and $8$ otherwise. 
\end{prop}
\Pf Every simple root has inner product at most $0$ with $n\rho$.
In a Lorentzian space the inner product of two vectors of nonpositive norm
in the same cone is at most $0$ and $0$ only if both vectors
are proportional to the same norm~$0$ vector. This implies that if we 
write $n\rho$ as sum of simple roots with positive coefficients the only 
simple roots appearing in the sum are positive multiples of~$\rho$.  
Since the support of an imaginary root
is connected it follows that all the~$n\rho$ are simple roots. 
Their multiplicities are given in Theorem~\ref{rootmult}. 
(Cf.~also Lemma~4 in section~3 of~\cite{B-fake}.)\eop

Now we show that we have already found all the simple roots of $\g$. 

\begin{thm}
A set of simple roots for $\g$ is the following. The real simple roots 
are the norm~$2$ vectors $\a$ in $L$ with $(\rho,\a)=-\a^2/2$ and the norm~$1$ 
vectors $\a$ in $L'$ with $(\rho,\a)=-\a^2/2$. The imaginary simple roots are
the positive multiples~$n\rho$ of $\rho$ with multiplicity~$16$ for even $n$  
and with multiplicity~$8$ for odd~$n$.
\end{thm}

\Pf The proof is analog to the proof of Theorem~7.2 in~\cite{B-moonshine}. 
Let ${\bf k}$ be the generalized Kac-Moody algebra with root lattice $L'$, 
Cartan subalgebra $L'\otimes {\bf R}$ and simple roots as stated in the theorem. 

In~\cite{S-twist}, Theorem~3.2, product and Fourier expansions of an automorphic 
form on the Grassmannian ${\rm Gr}_2(M \otimes {\bf R})$ with $M=L\oplus \II$ are 
worked out for different cusps by applying Borcherds' theory of theta lifts to the 
vector valued modular form $(f_{\gamma})_{\gamma\in M'/M}$. The expansion at the 
cusp corresponding to a primitive norm $0$ vector in the sublattice $\II\subset M$ 
shows that the denominator identity of ${\bf k}$ is given by
$$ e^{\rho}\prod_{\a\in L^+}\left(1-e^{\a}\right)^{c(-\a^2/2)}
\prod_{\a\in {L'}^+}\left(1-e^{\a}\right)^{c(-\a^2)}
\qquad\qquad\qquad\qquad\qquad $$
$$\qquad\qquad\qquad\qquad\qquad =
\sum_{w\in W} \det(w)\, w\left(e^{\rho}\prod_{n=1}^{\infty}
\big(1-e^{n\rho}\big)^8 \big(1-e^{2 n\rho}\big)^8\right).$$
Here, $W$ is the reflection group generated by norm~$1$ vectors of $L'$ and 
the norm~$2$ vectors of $L\subset L'$ and the exponents $c(n)$ are the 
coefficients of the modular form $h(\tau)=\sum_{n\in\Z}c(n)q^n$ defined 
in section~\ref{voav} (cf.~also~\cite{J}). 

Using Theorem~\ref{rootmult} and the definition of the $f_{\gamma}$,
we see that $\g$ and ${\bf k}$ have the same root multiplicities. 
The product in the denominator identity determines the simple roots of $\g$ 
because we have fixed a Cartan subalgebra and a 
fundamental Weyl chamber. It follows that $\g$ and ${\bf k}$ 
have the same simple roots and are isomorphic.
\eop

\begin{cor}
The denominator identity of $\g$ is 
$$ e^{\rho}\prod_{\a\in L^+}\left(1-e^{\a}\right)^{c(-\a^2/2)}
\prod_{\a\in {L'}^+}\left(1-e^{\a}\right)^{c(-\a^2)}
\qquad\qquad\qquad\qquad\qquad $$
$$\qquad\qquad\qquad\qquad\qquad =
\sum_{w\in W} \det(w)\, w\left(e^{\rho}\prod_{n=1}^{\infty}
\big(1-e^{n\rho}\big)^8 \big(1-e^{2 n\rho}\big)^8\right)$$
where $W$ is the reflection group generated by norm~$1$ vectors of $L'$ and 
the norm~$2$ vectors of $L$ and $c(n)$ is the coefficient of $q^n$ in
$$\eta(\tau)^{-8}\eta(2\tau)^{-8} = q^{-1}+  8+ 52\,q+ 256\, q^2 + 1122\, q^3 + 4352\, q^4 + \cdots .$$
\end{cor}

Using $L'(2)\cong\Lambda_{16} \oplus \II(2)$, we see that the denominator 
identity of $\g$ is a rescaled version of the denominator identity of 
Borcherds' fake baby monster Lie algebra determined
in~\cite{B-moonshine}, Sect.~14, Example~1. This implies:
\begin{cor}
The \gkm $\g$ is isomorphic to the fake baby monster Lie algebra.
\end{cor}

\medskip 

In a forthcoming paper we will describe similar constructions of some other 
generalized Kac-Moody algebras.


\end{document}